\documentclass[11pt]{article}
\usepackage[margin=1in]{geometry}
\usepackage{tgtermes}
\usepackage{setspace}  % For line spacing
\usepackage{subcaption}
\usepackage{changepage}
\usepackage{amsmath, amssymb, graphicx, natbib, url}
\usepackage{pgfplots}
\pgfplotsset{compat=1.18}
\usepgfplotslibrary{statistics}
\usepackage{pgfplotstable}
\usepackage[T1]{fontenc}
\usepackage[utf8]{inputenc}

\usepackage{microtype}

\usepackage[font=footnotesize, labelfont=bf, textfont=bf]{caption}
\usepackage{pgfplots}
\pgfplotsset{compat=1.18}
\usepackage{filecontents} 
\usepackage{algorithm}
\usepackage{algpseudocode}
\usepackage{tikz}
\usetikzlibrary{positioning, fit, backgrounds}
\usepackage[colorlinks=true,linkcolor=blue, citecolor= blue]{hyperref}%
\usepackage{natbib}
 \bibpunct[, ]{(}{)}{,}{a}{}{,}%

\usepackage{todonotes} 
\usepackage{xcolor}

\usepackage{threeparttable}  % In prea
\usepackage{fancyhdr} 
\usepackage{amsthm}

\newtheorem{theorem}{Theorem}     % Preferred (Theorem 1, Lemma 1, Theorem 2)
\newtheorem{definition}{Definition}[section]
\newtheorem{lemma}{Lemma}
\graphicspath{ {./info/}}
\title{Route Optimization Over Scheduled Services For Large-Scale Package Delivery Networks}

\author{Mohammed Faisal Ahmed \\ \small{Department of Industrial and Systems Engineering,
Georgia Institute of Technology, \texttt{mahmed91@gatech.edu}}\\Pascal Van Hentenryck \\ \small{Department of Industrial and Systems Engineering,
Georgia Institute of Technology, \texttt{pvh@gatech.edu}}\\Ahmed El Nashar \\ \small{United Parcel Service, \texttt{aelnashar1@ups.com}}}
\date{}
\begin{document}
\maketitle

% Enter all authors
%} % end of the block
\begin{abstract}
    This paper introduces the Trailer Path Optimization with
  Schedule Services Problem (TPOSSP) and proposes a column- generation
  heuristic (CG-heuristic) to find high-quality solutions to
  large-scale instances. The TPOSSP aims at determining trailer routes
  over a time-dependent network using existing scheduled services,
  while considering tractor capacity constraints and time windows for
  trailer pickups and deliveries. The objective is to minimize both
  the number of schedules used and the total miles traveled. To
  address the large scale of industrial instances, the paper proposes
  a network reduction technique that identifies the set of feasible
  schedule-legs for each requests.  Moreover, to address the resulting
  MIP models, that still contains hundred of millions variables, the
  paper proposes a stabilized column-generation, whose pricing problem
  is a time-dependent shortest path. The approach is evaluated on
  industrial instances both for tactical planning where
  requests for the entire network are re-optimized and for real-time operations where
  new requests are inserted.  In the tactical planning
  setting, the column-generation heuristic returns solutions with a
  3.7\%-5.7\% optimality gap (based on a MIP relaxation) in under 1.7
  - 10.3 hours, and improves the current practice by 2.3-3.2\%, with
  translates into savings of tens of millions of dollars a year. In
  the real-time setting, the column-generation heuristic returns
  solution within 3\% of optimality in under 1 minute, which makes it
  adequate for real-time deployment. The results also show that the
  network reduction decreases run times by 85\% for the
  column-generation heuristic.
\end{abstract}
\noindent \textbf{Keywords:} Package Delivery Network, Scheduled Services, Path Optimization, Column Generation, Lagrangian Dual, Network Reduction
%\HISTORY{Received: Month DD, YYYY; Accepted: Month DD, YYYY; Published Online: Month DD, YYYY}

%%%%%%%%%%%%%%%%%%%%%%%%%%%%%%%%%%%%%%%%%%%%%%%%%%%%%%%%%%%%%%%%%%%%%%

% Text of your paper here

\section{Introduction}\label{sec:Intro}

With the rapid growth of e-commerce, the volume of small shipments has
increased significantly, rising from 21.2 billions in 2022 to an
estimated 32 billions by 2028 (\cite{caponereport}). Small shipment
operations are costly and sensitive to economic fluctuations. Rising
fuel and labor costs drive up transportation expenses for Logistic
Service Providers (LSP), leading to higher prices for retail customers
and small businesses. Efficient resource management is therefore
critical to minimizing costs for both customers and companies.

LSPs that are operating in the small shipment business profit from
economies of scale. They can consolidate shipments from various
customers at hubs and move the shipments in trailers through a network
of these hubs. The trailers may be moved using different
transportation modes (e.g., trucks, rail, or air) depending on the
promised delivery due date of the shipment (also generally known as
service level). Transportation between two hubs in the network with or
without intermediate stops is referred to as a service. The LSP
planning process for timely and cost-effective delivery of shipments
is assisted by solving Service Network Design (SND) problems, first
introduced by \cite{crainic}. A SND problem involves the determination
of shipment paths, as well as the routes and schedules for the
transportation modes. Trailer movement frequencies are pre-planned
based on the expected shipment volume between origin-destination sort
pairs. A sort refers to the process of sorting shipments into the
correct trailers based on the shipment's next destination. Sorting
operations are carried out in dedicated sorting facilities during
specific shifts (e.g., night or day sorts). A planned trailer movement between sorting shifts of two facilities is called {\em
  a request} in this paper. A request consists of the following
attributes:
\begin{itemize}
    \item the origin hub where the trailer has to be picked up and
      the destination hub where the trailer has to be dropped;
    \item the earliest available time to pick up the trailer at the
      origin and latest arrival time for the trailer to arrive at the
      destination;
    \item The type of the trailer, e.g., a 53 foot trailer or a 28
      foot trailer.
\end{itemize}
Based on the requests, the driver schedules and fleet sizes are
optimized for the entire network.

A driver schedule typically comprises of a route starting and ending
at the driver's home hub, with specific start and end dates/times for
each leg of the route. Each leg of the route is referred to as a {\em
  schedule-leg} in this paper. Driver schedules must comply with
federal regulations regarding total work hours, meals, breaks, and
start and finish locations. Driver schedules are classified into two
types:
\begin{itemize}
    \item Same-day schedules: these operate within a single day, typically lasting 8 to 12 hours.
    \item Long-haul schedules: these may span up to a week and cover
      long-distance routes.
\end{itemize}

\begin{figure}[!t]
\centering
\includegraphics{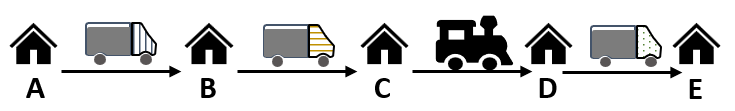}
\caption{Illustrating a Potential Path of a Trailer,}
\label{fig:1}
\end{figure}

Requests whose travel times are more than 6 hours cannot be fulfilled
by same-day schedules as it would exceed the daily max work hours for
their drivers. Such requests are either fulfilled just by long-haul
schedules, or by a combination of same-day schedules, or using both
same-day and long-haul schedules. Figure \ref{fig:1} illustrates a
path traversed by a request using multiple driver's schedule-legs. A
trailer begins at origin hub A and is transported to hub B on a
driver's schedule-leg. At hub B, the trailer is dropped and
temporarily parked until it is picked up by a different driver's
schedule-leg (indicated by a change in pattern of the tractor). The
trailer is then transported to rail yard C, moved by train to rail
yard D, and finally delivered to destination hub E by another driver's
schedule-leg. This path consists of one schedule-leg from three
different driver schedules and one rail schedule-leg.

A single tractor can transport multiple trailers
simultaneously. Consequently, trailers are often consolidated on
schedule-legs. Between schedule-legs, trailers are parked at exchange
points, where they wait until being picked up by another driver or
schedule, which allows for consolidation of trailers. The exchange
points can be sorting facilities, or dedicated exchange points. All
the exchange points are referred to as hubs in this paper.

\subsection{The TPOSSP and Its Motivation}
This paper introduces the Trailer Path Optimization using Scheduled
Services Problem (TPOSSP), which focuses on determining optimal paths
for requests by leveraging existing scheduled services. At a high
level, the TPOSSP can be defined as the problem of finding routes for
trailers using existing schedules, by considering tractor capacity
constraints on schedule-legs, and time constraints for earliest pickup
time and latest arrival time of requests. The route of a request
consists of sequence of schedule-legs where each subsequent leg must
depart after the arrival of the previous leg at a hub. The objective
of TPOSSP is to minimize the schedules used and the miles driven to
route all the requests.

{\em The TPOSSP is motivated by the need of LSPs to accommodate
  week-over-week variations in requests in a cost-effective manner
  without causing major disruptions to the scheduled operations}. Due
to fluctuations in shipment volume, there are weekly additions and
deletions of requests. Some changes are planned, while others are
ad-hoc and occur daily. In order to adapt to these changes, the LSPs
identify requests that regularly occur/do not occur, and adapt the
routes of trailer and schedules based on these request variations on a
weekly basis. {\em Planning and operations then face the challenge of
  finding ways to serve the requests, while keeping the changes in
  schedules at a minimum.} Keeping changes at a minimum stems from two
main reasons.
\begin{itemize}

\item LSPs are sometimes contractually obligated to provide consistent
  start times and routes for the drivers and frequent changes in
  routes negatively impact driver performance.
      
\item Schedules determine when the trailers arrive at the hub which
  affects the inbound trailer arrival time distribution into each
  sort. The  distribution of inbound trailer arrival times affects the
  personnel needed for unloading, sorting, and shifting in each
  sort. Therefore, schedule modifications affect the sort
  operations. Such changes may cause the sort operations to not
  function efficiently, which may result in delays at hubs and
  requests not meeting guaranteed service levels.

\item A complete rescheduling of the entire network for large LSPs on
  a weekly basis is a monumental task that would need a buy-in from
  all stakeholders, since it affects all aspects of the network.
\end{itemize}

\noindent
The resulting challenge is addressed with a two-stage approach. In the
first stage, the TPOSSP is solved to find paths for requests using
existing schedules. In the second stage, new schedules are created for
requests for which paths could not be found in the first stage. The
second stage problem can be formulated as a capacitated vehicle
routing problem with time windows, which is widely explored in the
literature. It is not the topic of this paper.

The TPOSSP itself can be solved at two levels: tactical
planning and real-time operations.
\begin{itemize}
\item The tactical planning TPOSSPs are large-scale and
  decide paths for all requests in the network.  The largest TPOSSP
  problems solved in this paper have number of requests multiplied by the number of schedule-legs given by \#requests x \#schedule-legs equal to 74 billion.

\item The real-time TPOSSPs are smaller scale problems concerned with
  finding paths for newly added requests, keeping the paths of
  existing requests fixed. The paths are found using the available
  capacity of legs after serving existing requests and removing
  canceled requests. The newly added requests have at their disposal
  the legs with available capacity throughout the network.  Real-time
  TPOSSPs contain fewer requests and have \#requests x \#schedule-legs in the range of 10 million to 80 million. To meet the real-time requirements of
  operations, they must be solved with 1 minute or so.
\end{itemize}

\subsection{Methodological Approach}

To address the large scale of industrial TPOSSP instances, the paper
first proposes a network reduction technique that identifies the set
of feasible schedule-legs for each request.  Moreover, to address the
complexity of the resulting MIP models, that still contain hundred of
millions variables, the paper proposes a stabilized column-generation,
whose pricing problem is a time-dependent shortest path. The approach
is evaluated both in tactical planning where the entire network is re-optimized and in a real-time setting where only a few fraction of requests are
inserted.  Experimental results on industrial instances show that the
network reduction decreases run times by 53\%–85\% for the
column-generation heuristic. In tactical planning setting, the
column-generation heuristic returns solutions with a 3.7-5.7\%
optimality gap (based on a MIP relaxation) in under 1.7-10.3 hours,
and improves the current practice by 2.3-3.2\%. In real-time setting,
the column-generation heuristic returns solution within 3\% of
optimality in under 1 minute, which makes it adequate for real-time
deployment.

\subsection{Contributions} The main contributions of the paper can be summarized as follows:
\begin{enumerate}
\item The paper formalizes the TPOSSP, a critical aspect of LSP
  operations, for the first time.
\item The paper proposes a network reduction method that significantly
  reduces the size of the problem, thereby decreasing the runtimes of
  both the exact arc-based formulations by 85\%-99\% and the column
  generation heuristic by 53\%-85\%.
\item The paper proposes a stabilized CG-heuristic whose pricing
  problem is a time-dependent shortest path algorithm for solving
  large-scale TPOSSP.
\item Experiments on real world tactical and real-time instances
  demonstrate the scalability and performance of the CG-heuristic. The
  large-scale data consisted of \#requests x \#schedule-legs nearly equal to 74 billion. They also demonstrate that the
  CG-heuristic meets the run-time requirements for real-time
  operations.
\item The experiments show that the CG-heuristic solutions reach
  within 3.7\%-5.7\% of the true optimal solution for tactical
  TPOSSPs, a 2.3\%-3.2\% reduction in total operating cost, a
  8.5\%-10.3\% reduction in empty miles (miles driven without a
  trailer). This translates to {\em millions of dollars in net cost
    reduction per week}.  On real-time instances, the CG-heuristic
  typically produces solutions with 3\% of optimality under 1 minute.
\end{enumerate}

The rest of this paper is organized as follows. Sections
\ref{sec:LitRev}, \ref{sec:Pdef}, and \ref{section:arc-formulation}
present the literature review, the problem definition, and the
arc-based formulation. Section \ref{section:reduction} presents the
network reduction technique, and Section \ref{section:path} the
path-based formulation and column generation heuristic. Section
\ref{section:path} describes the computational experiments and Section
\ref{section:conclusion} concludes the paper. Proofs are given in the
Appendix.

\section{Literature Review}\label{sec:LitRev}

Extensive literature searches in operations research, transportation,
and logistics revealed that the specific TPOSSP has not
been addressed. Problems that have similar mathematical models to
the TPOSSP have been explored in the SND literature, specifically in rail
network optimization. This section reviews related literature in
Service Network Design (SND), Multi-Commodity Network Flow (MCNF), and
Time-Dependent-Shortest-Path-Problem (TDSPP).

SND was first introduced by \cite{crainic} and
\cite{Wieberneit}. \cite{crainic} framed the service network design
problem as an MCNF problem, providing both arc and path-based
formulations. \cite{cranic2} solved a variant of the SND problem that
jointly solved for paths and schedules which is related to TPOSSP but
varies in the fact that TPOSSP does not solve for schedules. They
modeled the problem on a time-expanded network and developed a column
generation-based meta-heuristic with two steps: first, solving the
linear relaxation using column generation, which involves a pricing
problem to generate least-cost cycles for scheduling resources;
second, using these cycles to find good integer feasible solutions
through a slope scaling procedure (\cite{kim}). Compared to
\cite{cranic2}, this paper considers instances that are an order of
magnitude greater, avoids the network discretization, and includes an
optimal network reduction technique to meet the requirements of
industrial instances.

\cite{zhu} solved tactical planning for ``freight rail carriers aiming
to select the train services to operate over the contemplated schedule
length (e.g., the week), together with their frequencies or schedules
(timetables), the blocks that will make up each train, the blocks to
be built in each terminal, and the routing of the cars loaded with the
customer's freight using these services and blocks'' (movements are
also considered). They presented an integrated schedule network design
methodology for rail freight transportation. Their proposed framework
contains a 3-layer space-time network of the associated operations and
decisions made. They discretized time and added hold-over arcs at the
yards to model time-dependent decisions. They proposed a
meta-heuristic name ellipsoidal search that explores large
neighborhoods of good solutions using information from the history of
search. The TPOSSP can also be modeled using hold-over arcs at the
hubs to make the network time independent. Adding hold-over arcs at
the hubs however would triple the number of legs (one hold-over arc
for each incoming leg at a hub and one for the outgoing leg). Even for
real-time problems, tripling the legs would blow up the number of legs
to millions. This paper proposes methods for pruning and quickly
solving the sub-problems on time-dependent networks, so there is no
benefit in discretizing the problem.

The TPOSSP can be more generally formulated as a time-dependent
multi-commodity network flow (MCNF) problem with schedule
constraints. The MCNF problem has been used to model many
transportation and logistics problems (e.g., \cite{ayar},
\cite{Lienkamp}). The most important exact methods for solving MCNFs
are branch-and-cut and decomposition methods. Branch-and-price, a
popular decomposition approach developed by \cite{cynthia2}, was used
to solve multi-origin-destination integer multi-commodity flow
problems by \cite{cynthia}. Popular heuristic methods for solving MCNFs
include simulated annealing, tabu search, and genetic algorithms
(\cite{Salimifard2022}).
 
Some studies related to the problem in this paper include those by
\cite{ayar} and \cite{Lienkamp}. \cite{ayar} explored a
multi-commodity flow problem where commodities must be picked up and
delivered within a time window using trucks and ships, aiming at
minimizing transportation and stocking costs. They proved that the
problem is NP-hard and developed strong formulations. \cite{Lienkamp}
modeled passenger transportation on scheduled bus and train services
in an urban setting as an MCNF, aiming at minimizing transit
times. They used column generation to solve the linear relaxation,
branch-and-price for integer feasible solutions, and the $A^*$
algorithm for the pricing problem.

Most studies on time-dependent multi-commodity network flow problems
transform the network into a static one by creating a time-expanded
graph (e.g., \cite{cranic2}, \cite{Lienkamp}, \cite{ayar}). For large-scale
networks, this transformation can lead to much larger
networks. \cite{boland} introduced partially time-expanded networks
for SND, where not all time steps are included, allowing for iterative
refinement to reduce graph size.

In this paper, the pricing problem of column generation reduces to a
TDSPP with step-wise arc cost functions and waiting allowed at
nodes. \cite{he} introduced a dynamic discretionary discovery
algorithm to solve TDSPPs with piecewise linear arc travel time
functions. This algorithm operates within partially time-expanded
networks, with arc costs representing the minimum travel time over the
next interval. The shortest path in this network offers a lower bound
on the optimal path’s value, while upper bounds are obtained during
lower bound calculations. The algorithm iteratively improves
discretization using breakpoints in the arc travel time functions,
refining both lower and upper bounds and ultimately proving optimality
after a finite number of iterations.

Given the large network sizes explored in this paper, it is beneficial
to exclude legs that will never be visited by a feasible
path. Literature on On-Demand Multi-Modal Transit Systems (ODMTS)
introduced by \cite{maheo} has used pre-processing techniques to
eliminate arcs that will not be used by feasible paths. More recent
research by \cite{guan} for ODMTS design incorporated a sophisticated
approach to eliminate hubs that will not lie on an optimal path by
finding shortest paths from the origin of a request to the hub and
from the hub to the destination. If the cost of the path exceeds the
upper bound of the optimal cost, the hub is eliminated.

Based on this review, it can be concluded that, from an application
point of view, the TPOSSP is a new problem, not yet explored in
current literature, although similar optimization models exist in train
network optimization problems. The TPOSSP can be formulated as time
dependent MCNF problem, and solution approaches used for MCNF are also
applicable to the TPOSSP. However, the methodology proposed in this
paper to solve TPOSSPs is unique, motivated by its specific
application and the large scale nature of its instances. In
particular, this paper does not discretize the network, proposes an
algorithm to find the optimal sub-network relevant to each request,
and leverages a stabilized column generation with a time-dependent
pricing subproblem.

\section{Problem Definition} \label{sec:Pdef}

\begin{figure}[!t]
\centering
\includegraphics{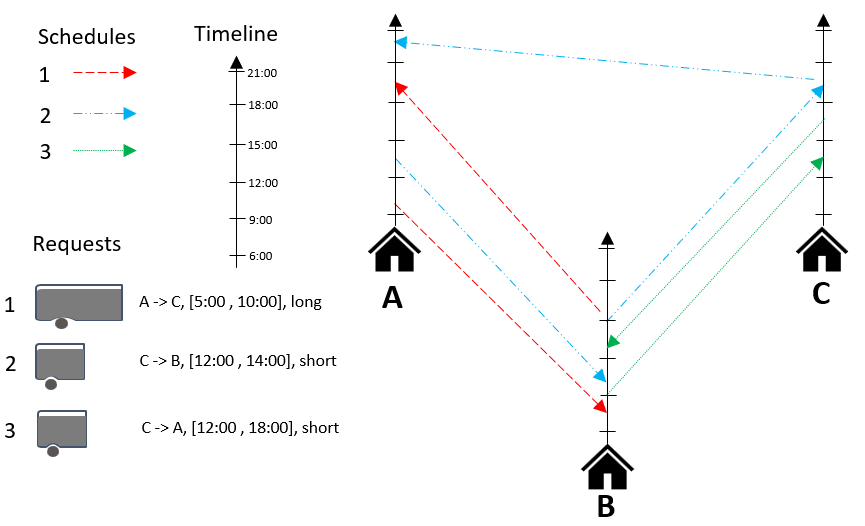}
\caption{A TPOSSP Instance with three schedules and three inputs and three hubs. Each schedule is shown with colored arrows and the start-end location, earliest available and latest arrival times, and the type of trailer is shown for each request. }
\label{fig:2}
\end{figure}

\begin{figure}[!t]
\centering
\includegraphics{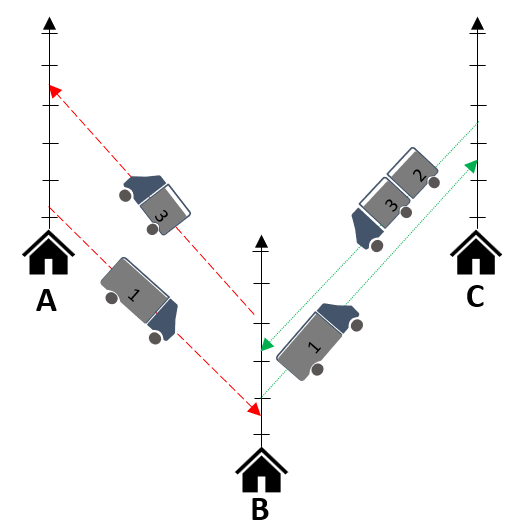}
\caption{A Solution to the TPOSSP Instance: The solution contains the path of each request, and schedule 2 was eliminated.}
\label{fig:3}
 \end{figure}

This section formally introduces the TPOSSP problem. To build
intuition, Figure \ref{fig:2} presents an illustrative instance of the
TPOSSP problem. There are three schedules and three requests. The
origin, destination, time windows, and trailer type are shown next to
the requests. A solution is shown in Figure \ref{fig:3}. Request 1
uses the first schedule-leg of Schedule 1 and the first schedule-leg
of Schedule 3. Requests 2 and 3 use short trailers and can be pulled
by a single tractor, so they are combined in the second leg of
Schedule 3. Request 3 completes its journey by utilizing the second
leg of Schedule 1. None of the requests use Schedule 2, so it was
eliminated thereby minimizing the cost.

\begin{table}[!t]
\begin{threeparttable}
\small
\caption{The Notations Used in the Paper.\label{tab:nomenclature}}
{\begin{tabular}{@{}l@{\quad}l@{}}
\hline 
Symbol              & Definition           \\ \hline
Sets: &  \\
     & \\
    $R$ & Set of loaded trailer requests\\
    $L$ & Set of schedule-legs  \\ 
    $S$ & Set of schedules  \\ 
    $H$ &  Set of hubs\\
    $H^h_{in}$ & Set of incoming legs of hub $h$\\
    $H^h_{out}$ & Set of outgoing legs of hub $h$\\
    $H^r_{inter}$ & Set of hubs other than origin and destination hub of request $r$ \\
    $O^r$ & Set of origin schedule-legs of request $r$\\
    $D^r$ & \textrm{Set of destination schedule-legs of request $r$}\\
    $K^r$ & \textrm{Set of schedule-legs except origin and destination schedule-legs of request $r$ }\\
     & \\
    Parameters: &\\
     & \\
    $\sigma_{s}$ & Fixed cost of cutting schedule $s$\\
    $\theta_l$ & Cost per mile per short trailer for leg $l$\\
    $d_l$ & Start hub of leg $l$\\
    $e_l$ & End hub of leg $l$\\
    $t^{s}_l$ & Start time of leg $l$\\
    $t^{e}_l$ & End time of leg $l$ \\
    $s_l$ & Schedule of leg $l$ \\
    $c_l$ & Capacity of leg $l$ \\
    $m_l$ & Miles of leg $l$\\
    $p_r$ & Origin hub of request $r$\\
    $q_r$ & Destination hub of request $r$\\
    $t^{est}_r$ & Earliest start time of request $r$ from the origin\\
    $t^{lat}_r$ & Latest arrival time of request $r$ to the destination\\
    $v_r$ & Short trailer equivalents of request $r$\\ \\ \hline

\end{tabular}}
\begin{tablenotes}
\footnotesize
\item Variables are defined with the formulations.
\end{tablenotes}
\end{threeparttable}
\end{table}

The goal of the TPOSSP is to find paths for each request using
pre-defined schedules in such a way that the total number of schedules
utilized and the total miles driven are minimized. The notations used
throughout the paper are shown in Table \ref{tab:nomenclature}. The
inputs to TPOSSP are the set of requests $R$, schedules $S$,
schedule-legs (also referred to as legs) $L$, and hubs $H$. Each
request $r$ has an origin $p_r$, destination $q_r$, and a time window
defined by the earliest allowed start time at the origin $t^{est}_r$
and latest allowed arrival time at the destination $t^{lat}_r$. Each
leg $l$ of a schedule has an origin $p^s_l$, a destination $p^e_l$, a
start time from the origin $t^s_l$, and an arrival time at the
destination $t^e_l$. The set of incoming and outgoing legs of hub $h$
are denoted by $H^h_{in}$ and $H^h_{out}$, respectively. The set
$H^r_{inter}$ includes hubs other than the origin and destination of
request $r$. The set of outbound legs from the origin of request $r$
is denoted by $O^r$, the set of inbound legs to the destination of
request $r$ is given by $D^r$, and the set of legs other than those
$O^r$ and $D^r$ is denoted by $K^r$. The travel time on every
schedule-leg is strictly positive.
 
The algorithms and models developed in this paper are applicable to
any mode of transportation, but the paper focuses on trucks. Each
request corresponds to a single trailer, which can vary in dimensions,
e.g., 53' (long), 45', 48', and 28' (short). A tractor can carry
either one long trailer or up to three short trailers. It is also
possible to combine two 45' trailers, or one 28' and one 45' trailer,
among other combinations. The volume of a short trailer (28') is
normalized to 1, and the volumes of other trailer types are defined
based on their equivalent short trailers. The capacity of a tractor
can be up to 3, meaning that it can carry up to three short
trailers. The capacity of a leg $c_l$ is the number of short trailers
that can be carried on that leg. The volume of a request $v_r$ is also
given in short trailer equivalents. For example, the volume is
normalized to 2.5 for a long trailer (which cannot be combined with
any other trailer). Intermediate trailer types such as 45' are
normalized to 1.5, and 48' is normalized to 1.9, as it is possible to
combine two 45' trailers, one 48', and one 28' trailer, but not two 48'
trailers.

The total cost is divided into two parts: the schedule cost and the
mile cost. The cost of schedule $s$ is given by $\sigma_s$, which
includes the driver cost and the mileage cost of moving the tractor
(without trailers) on the route. The mile cost is the additional cost
of moving a trailer on a leg. The cost of moving a short trailer over
a unit distance is given by $\theta_l$. More precisely, the total mile
cost of moving a request $r$ on a leg $l$ would be $\theta_l*m_l*v_r$.
It is assumed that the start and end times of legs $t^s_l$ and $t^e_l$
include the time required to attach and detach trailers to and form
the tractors at the hubs.
 
Requests that have no feasible path in a given schedule are removed:
they will need their own dedicated schedules. For the remaining
requests, the problem specification adds direct schedules from the
origin to the destination of every request and back in order to ensure
the model is always feasible (since not all requests may always be
served by the given schedule). The schedule cost of the added direct
schedules is the cost of creating a new schedule from the origin of
the load to the destination and back, while the mile cost of their
schedule-legs is the cost of moving directly from the origin of the
load to its destination and vice versa. The cost of an added schedule
is given by $\sigma^{dummy}_s$.  The objective function of the TPOSSP
is to minimize the schedule and mile cost.

The constraints of the TPOSSP can be categorized into three categories.

\begin{itemize}

\item \textbf{Path feasibility}: These constraints make sure that
   there is at least one path for every request and the path is time
   feasible.
   
\item \textbf{Capacity constraints}: These constraints make sure
  that each leg capacity is not exceeded.

\item \textbf{Schedule identification constraints}: These constraints
  identify which schedules are being utilized.
 \end{itemize}

\section{Arc Based Formulation}
\label{section:arc-formulation}

The TPOSSP can be formulated as a Mixed Integer Program (MIP). There
are two types of decision variables for the MIP.  Variable $x_{lr}$
(\textbf{x} in vector format) represents whether request $r$ uses
schedule-leg $l$ in its path. Variable $y_s$ (\textbf{y} in vector
format) represents whether schedule $s$ is used.

\begin{figure}[!t]
    \begin{subequations}
\begin{align}
   &\sum_{l \in H^h_{in}}x_{lr}= 0 \quad  \forall r \in R ,h = p_r \label{g1}\\
   &\sum_{l \in H^h_{out}}x_{lr}= 0 \quad  \forall r \in R ,h = q_r \label{g2}\\
   &\sum_{l \in H^h_{in}}x_{lr}- \sum_{l \in H^h_{out}}x_{lr} =\begin{cases}
            -1, & \text{if $h = p_r$ , $r \in R$ } \\ 
			0, & \text{if $h \in H^r_{inter}$ , $r \in R$ }\\
            1, & \text{if $h = q_r$ , $r \in R$}
		 \end{cases} \label{g3}\\
   &  \sum_{l \in H^h_{in}} t^{e}_l * x_{lr} \leq \sum_{l \in H^h_{out}} t^{s}_l * x_{lr} \quad \forall h \in H^r_{inter}, \forall r \in R \label{g7}  \\
   & t^{est}_r \leq  \sum_{l \in L^r_{origin}} t^{s}_l* x_{lr} \quad \forall r \in R  \label{g8} \\
   & \sum_{l \in L^r_{destination}} t^{e}_l* x_{lr} \leq t^{lat}_r \quad \forall r \in R \label{g9} 
\end{align}
\end{subequations}
\caption{The Path Feasibility Constraints $\mathcal{G}(\textbf{x})$.} \label{pathfeasibilityconstr}
\end{figure}

The path feasibility constraint set is denoted by
$\mathcal{G}(\textbf{x})$ and given in Figure
\ref{pathfeasibilityconstr}. Constraints (\ref{g1}) and (\ref{g2})
ensure that there are no incoming legs to the origin of a request and
no outgoing legs from the destination of a request. Constraints
(\ref{g3}) ensure that there is at least one outgoing leg from the
origin hub and at least one incoming leg to the destination hub for
each request; it also captures the flow balance for all intermediate
hubs. Constraint (\ref{g7}) ensures that the incoming leg to a hub arrives
before the outgoing leg departs. Constraints (\ref{g8}) and (\ref{g9})
ensure that the path satisfies the time window constraints of the request.

The capacity constraints ensure that the tractor capacity is not
exceeded. The capacity constraint set is denoted by
$\mathcal{C}(\textbf{x},\textbf{y})$ and includes the following
constraints:
\begin{equation}
    \sum_{r \in R} v_r * x_{lr} \leq c_{l}*y_s \quad \forall l \in L, s= s_l \label{c1}
\end{equation}
The schedule constraints ensure that a schedule is selected if at
least one of its schedule-leg is used. The set of schedule constraints
is denoted by $\mathcal{E}(\textbf{x},\textbf{y})$ and contains the
following constraints:
\begin{equation}
    x_{lr} \leq y_s \quad \forall l \in L, \forall r \in R, s= s_l \label{s1}
\end{equation}

\begin{figure}[!]
\setlength{\abovedisplayskip}{0pt}
\setlength{\belowdisplayskip}{0pt}
   \begin{subequations}
\begin{align*}
    &\min_{\textbf{x}, \textbf{y}} \quad \sum_{s \in S} \sigma_s*y_s+ \sum_{l \in L} \sum_{r \in R} \theta_l*m_l*v_r*x_{lr}
\end{align*}
\begin{equation}
\text{s.t. } \{\textbf{x},\textbf{y}\} \in  \mathcal{G(\textbf{x})}\cap \mathcal{C(\textbf{x},\textbf{y})} \cap  \mathcal{E(\textbf{x},\textbf{y})} \label{arc1}
\end{equation}
\end{subequations}
\begin{align*}
&  \textbf{x} \in \{0,1\}^{|L|X|R|}, \textbf{y} \in \{0,1\}^{|S|}
\end{align*}
\caption{The Arc-Based MIP Formulation for the TPOSSP
  (Arc-MIP).} \label{arcmip}
\end{figure}
The schedule constraint set $\mathcal{E}$ is redundant for the above
problem because these constraints are already captured by
$\mathcal{C}$. However, $\mathcal{E}$ is included because it is not
redundant in the linear relaxation, thereby providing a stronger lower
bound, which is beneficial computationally.  The full arc-based
formulation is given in Figure \ref{arcmip}.

\section{Optimal Sub-Networks}
\label{section:reduction}

The TPOSSP is solved at the global network level. For the industry
cases motivating this paper, the instances have \#requests x \#schedule-legs up to 80 million,
even for real-time operations. To make the problems tractable, it is
beneficial to reduce the network size to include only these legs that
lie on some feasible path for each request. This section proposes an
exact polynomial time algorithm to identify these legs and prune the
rest of the network. The experimental results show that leg pruning
can reduce the run times of the CG-heuristic by 53\%-85\%.

\begin{definition}[Optimal Sub-Network]
An optimal sub-network for request $r$ is the union of the legs that
appear on a feasible path for $r$, i.e., a path satisfying the
path-feasibility constraints $\mathcal{G}$.
\end{definition}

This section presents an exact polynomial-time algorithm, known as the
Earliest and Latest Start Time Pruning (EALSP) algorithm, to identify
the optimal sub-network of a request. For request $r$, the EALSP
algorithm works on a preprocessed graph $G^r = (H^r,L^r)$ obtained
from the original graph $G=(H,L)$ by the following transformations:
\begin{itemize}
    \item Remove legs $l \in L$ that do not have enough capacity, i.e., $v_r > c_l$.
    \item Remove the incoming legs to the origin of request $r$ and
      the outgoing legs from the destination of request $r$.
    \item Remove legs $l \in L$ where $t^s_l < t^{est}_r$ or $t^e_l > t^{lat}_r$.
\end{itemize}

\begin{table}[!t]
\small
\begin{adjustwidth}{-1.5cm}{0cm}
\centering
\caption{Notations Used in Sub-network\label{tab:nomenclaturesubnet}}
{\begin{tabular}{@{}l@{\quad}l@{}}
\hline
Symbol              & Definition           \\ \hline
     & \\
    $E^r_{*}$ & \textrm{The set of legs in the optimal sub-network of request $r$}.\\
    $H^r_{*}$ & \textrm{ Origin and destination hubs of legs in $E^r_{*}$}.\\
    $G^r_{*}$ & \textrm{ Optimal sub-network of request $r$ = ($H^r_{*}$,$E^r_{*}$)}.\\
    $eat^r_i$ & \textrm{For a request $r$, $eat^r_i$ denotes the earliest arrival time of a path departing from origin of request $r$ after $t^{est}_r$ and arriving at hub $i$}.\\
    &\textrm{If there is no feasible path from origin of request $r$ to hub $i$ then $eat^r_i = \infty$. The vector form for a request $r$ is denoted as $\mathbf{eat^r}$}.\\
    $lst^r_i$ & \textrm{For a request $r$, $lst^r_i$ denotes the latest start time of a path departing from hub $i$ and arriving at destination of request $r$ before $t^{lat}_r$.}\\
    & \textrm{If there is no feasible path from hub $i$ to destination of request $r$, then $lst^r_i = -\infty$. The vector form for a request $r$ is denoted as $\mathbf{ lst^r}$}.\\
    $E^r$ & \textrm{The set of legs whose start time is less than $eat^r_i$ and end time is less than $lst^r_j$.}\\
    & \textrm{Where $i$ is the origin of the leg and $j$ is the destination of leg.}\\
    $H^r$ & \textrm{Origin and destination hubs of legs in $E^r$}\\
\hline
\end{tabular}}{}
\end{adjustwidth}
\end{table}
\subsection{Sub-Network Properties}

The idea behind the EALSP algorithm is to find the earliest arrival
date of request $r$ at hub $i$, denoted by $eat^r_i$, and the latest
starting time of request $r$ from hub $i$, denoted by $lst^r_i$.
The optimal sub-network for request $r$ can then be defined as

\setlength{\abovedisplayskip}{0pt}
\setlength{\belowdisplayskip}{0pt}
\begin{subequations} \label{nreq}
\begin{align}
  E^r = \{l \in L: eat^r_{d_l} \leq t^s_{l}, t^e_{l} \leq lst^r_{e_l} \}.
\end{align}
\end{subequations}

\noindent
Note that $E^r = \emptyset$ if there is no feasible path from origin
to destination of request $r$.

\begin{lemma} \label{feaspathlemma}
Assume that a leg $l$ with start hub $i$ and end hub $j$ satisfies the
condition $eat^r_i \leq t^s_l \leq lst^r_j$. Then $l$ lies on a
feasible path from the origin to the destination of request $r$.
\end{lemma}
\begin{theorem}
  \label{thm1}
The optimal sub-network for request $r$ is $G^r_* = (H^r,E^r)$.
\end{theorem}

%The implication of Theorem \ref{thm1} is that, if $eat^r_i$ and $lst^r_i$ is found for every hub in the network, then the sub-network can be constructed using Equation (\ref{nreq}). The next section shows an algorithm to find $eat^r_i$ and $lst^r_i$ and eventually find the sub-network.\\

\subsection{Earliest Arrival Time and Latest Start Time Pruning Algorithm (EALSP)}

The EALSP algorithm finds the optimal sub-network $G^r_{*}$ for each
request $r \in R$ independently and in parallel.  It contains three
stages:
\begin{itemize}
    \item Find ${\bf lst}^r$ using Latest Start Time Pruning (LSP) algorithm.
    \item Find ${\bf eat}^r$ using Earliest Arrival Time Pruning (EAP) algorithm.
    \item Find $E^r_*$ using Equation (\ref{nreq}).
\end{itemize}
The LSP algorithm uses a dynamic programming approach. It performs the initialization
\setlength{\abovedisplayskip}{0pt}
\setlength{\belowdisplayskip}{0pt}
\begin{align*}
& lst^r_{q_r} = t^{lat}_r \\
& lst^r_{i} = -\infty \;\; (i \neq q_r),
\end{align*}
and applies the update rule 
\begin{equation*} 
lst^r_{h} = \max_{l \in H^h_{out}} \{ t^s_l | t^e_l \leq lst^r_{e_l} \} 
\end{equation*}
for each hub $h$ until no update takes place. In practice, the LSP
algorithm uses a Dijkstra-like strategy, starting at the destination
and exploring the hub backwards by maintaining the set of explored and
unexplored hubs, organized as a priority queue.  This process
continues iteratively until either the origin $p^r$ of request $r$ is
reached or until $lst^r_h < t^{est}_r$. If hub $p^r$ is not reached by
the end of the process, no feasible path exists.

The EAP algorithm is a mirror of LSP. It performs the initialization
\setlength{\abovedisplayskip}{0pt}
\setlength{\belowdisplayskip}{0pt}
  \begin{align*}
  &eat^r_{p_r} = t^{est}_r\\
  &eat^r_{i} = \infty \;\; (i \neq p_r),
\end{align*}
and applies the update rule 
\begin{equation*}
eat^r_{h} = \min_{l \in H^h_{in}} \{ t^e_l | t^s_l \geq eat^r_{d_l} \} \label{recbest}.
\end{equation*}

\section{Path-Based Formulation and Column Generation}
\label{section:path}

The largest instances motivating this study involves \#requests x \#schedule-legs nearly equal to 74 billion, which is also the number of variables in the problem. Even after reducing the
network size using the optimal sub-network, the arc-based formulation
still contains nearly 160 million variables for the large-scale
instances, which motivates the use of decomposition techniques. This
section presents a path-based formulation of the TPOSSP that uses a
Dantzig-Wolfe decomposition (\cite{Dantzig}) with $|R|$
subproblems. It also proposes a polynomial-time algorithm to solve the
pricing subproblem and a column generation-based heuristic to produce
high-quality solutions in reasonable time. The quality of the column
generation is quantified using a Lagrangian lower bound.

\subsection{Path-Based Formulation}

The Dantzig-Wolfe decomposition keeps the constraint sets
$\mathcal{C}$ and $\mathcal{E}$ in the master problem, and decomposes
the constraint set $\mathcal{G}$ per request. As is traditional, a
restricted master problem with only a subset of paths is solved. New
paths are added to the restricted problem by solving a pricing problem
that finds a feasible path with a negative reduced cost. The
formulation of the restricted master problem is shown in Figure
\ref{rmpformulation}. Here, $x^r_{pl}$ is an input parameter
indicating whether leg $l$ is used by path $p$ of request $r$,
$\lambda^r_p$ is the multiplier of path $p$ for request $r$, and
$Cost^r_p$ is the mile cost of path $p$.

\begin{figure}[!t]
\setlength{\abovedisplayskip}{0pt}
\setlength{\belowdisplayskip}{0pt}
 \begin{align*}
     \min_{  \boldsymbol{\lambda}, \textbf {Y}} \sum_{s \in S} \sigma_s*y_s+ \sum_{p \in P^r} Cost^r_p \lambda^r_p 
\end{align*}
\begin{subequations}\label{eq:litdiff}
\begin{align}
    \text{s.t.} &\sum_{r \in R} \sum_{p \in P} v_r x^r_{pl} \lambda^r_p \leq c_l*y_s \quad \forall l \in L , s=s_l \quad [\pi^c_l]  \label{rmp1} \\
     &\sum_{p \in P} x^r_{pl} \lambda^r_p \leq y_s \quad \forall r \in R, \forall l \in L , s=s_l \quad [\pi^s_{lr}]  \label{rmp2} \\
     &\sum_{p \in P} \lambda^r_p = 1 \quad \forall r \in R \label{rmp3} \quad [\pi^r]
\end{align}
\end{subequations}
\begin{align*}
\textbf{y} \in [0, 1]^{|S|},\boldsymbol{\lambda} \in [0, \infty)^{|R|X|P|}
\end{align*}
\caption{The Restricted Master Problem using Path Sets $P^r$ for
  Request $r$. The dual values of each constraint are given in
  brackets near the constraint.} \label{rmpformulation}
\end{figure}

\begin{figure}[!t]
\setlength{\abovedisplayskip}{0pt}
\setlength{\belowdisplayskip}{0pt}
   \begin{align*}
     \min_{\textbf{x} } \sum_{l \in L} \theta_l*m_l*v_r*x_{lr} - \\ \sum_{l \in L} (\pi^c_l*v_r+\pi^s_{lr})*x_{lr}-\pi^r
    \end{align*}
    \begin{align}
     \text{s.t. : }&\textbf{x} \in  \mathcal{G}(x),\textbf{ x} \in  \{0,1\}^{|L|}
    \end{align}
    \caption{The Sub-Problem Formulation.} \label{subprobformulation}
\end{figure}

The subproblem formulation is shown in Figure
\ref{subprobformulation}. Variable $x_{lr}$ indicates whether leg $l$
is used by the path of request $r$. The subproblem can be formulated
as a Time-Dependent-Shortest-Path-Problem (TDSPP) with non-negative
leg cost, where the time-dependency comes from the fact that the
possible outbound leg options from a hub depend on the arrival time to
the hub. The cost of the leg is given by $\alpha_{lr} =
\theta_l*m_l*v_r - (\pi^c_l*v_r+\pi^s_{lr})$. Note that $\alpha_{lr}
>0$, since $\theta_l*m_l*v_r > 0$, $\pi^c_l< 0$, and $\pi^s_{lr} < 0$.
The TDSPP is again solved by a Dijkstra-like algorithm. For every leg
$l$ of request $r$, the TDSPP algorithm computes the cost $\zeta_{lr}$
of the shortest feasible path starting from $l$ and ending at the
destination of request $r$. It performs the initialization
\setlength{\abovedisplayskip}{0pt}
\setlength{\belowdisplayskip}{0pt}
\begin{align*}
  \zeta_{lr} &= \alpha_{lr}, \;\; (l \in H^{q_r}_{in})\\
\zeta_{lr} &= \infty \;\; (l \notin H^{q_r}_{in}) 
\end{align*}
and applies the update
\begin{equation*}
\zeta_{lr} = \alpha_{lr}+\min_{l^{'} \in H^{e_l}_{out}}\{ \zeta_{l^{'}r} | t^e_l \leq t^s_{l^{'}} \} \;\;\ (l \in L \setminus H^{q_r}_{in}).
\end{equation*}
The algorithm can be enhanced with leg pruning based on arrival time
cost function. The arrival time cost function $f^r_h(t)$ of hub $h$
and request $r$ measures the cost of the shortest feasible path
starting from $h$ after time $t$ and ending at $q_r$, i.e.,
\begin{align*}
    f^r_h(t) = \min_{l \in H^h_{out}}(\zeta_{lr} | t \leq t^s_l \text{ \& } d_l = h) 
\end{align*}
and is step-wise left-continuous non-decreasing. The times $t$ at
which the $f^r_h(t)$ jumps (i.e., changes values) are called {\em
  critical points} and the legs departing at critical points are
called {\em critical legs}.  Lemma \ref{lemma2} implies that legs that
do not depart on critical points need not be explored by the TDSPP as they
are not essential to find the shortest path.
\begin{lemma}\label{lemma2}
There always exists an optimal path for the pricing subproblem that uses only the critical legs.
\end{lemma}

\subsection{Column Generation Stabilization} \label{cgstablization-section}

In large-scale column generation problems, dual values may not
converge smoothly and may oscillate frequently (\cite{Lübbecke}). This
paper uses Weighted Dantzig-Wolfe (WDW) decomposition by
\cite{Wentges} to stabilize the column generation. Its key idea is to
``“search for good Lagrangian multipliers in the neighborhood of the
best multipliers found so far'' (\cite{Wentges}). Instead of solving
the pricing subproblem with the optimal dual variables $\pi$, WDW
adjusts them through a convex combination with the best dual values
found so far. The adjusted duals $\Bar{\pi}^{k+1}$ at step $k+1$
are computed using 
\begin{align*}
\Bar{\pi}^{k+1} := \frac{1}{w_k}{\pi}^{k+1}+\frac{w_k-1}{w_k}\Bar{\pi}^{best,k} 
\end{align*}
where $\Bar{\pi}^{best,k}$ is the adjusted dual values of the best
objective lower bound found so far, and $w_k$ is a weight based on the number of
latest consecutive dual objective improvements. The columns found by
the pricing sub-problem using the adjusted values are then evaluated
using the actual dual values, and columns that have non-negative
reduced cost using the original dual values are eliminated. If none of
the columns have negative reduced cost, the algorithm shifts to the
regular Dantzig-Wolfe decomposition to converge to the true optimal
value.

%\begin{multline} 
%\label{eq:weightupdate}
%w_k := min(const, (k+ improvs)/2)
%\end{multline}
%\end{subequations}

\subsection{Finding Integer Solutions}

To find integer solutions, the CG-heuristic solves the restricted
master problem with the paths available at the convergence of the
column generation and the the integrality constraints on ${\bf y}$
variables.

\subsection{The Lagrangian Dual Bound}

This paper also uses a Lagrangian dual formulation to obtain strong
lower bounds. The two coupling constraints in the arc based
formulation $\mathcal{C}$ and $\mathcal{E}$ are dualized and their
Lagrangian multipliers are denoted by $\mathcal{C}$ by
$\boldsymbol{\kappa^C}$ and $\boldsymbol{\kappa^E}$ respectively. The
Lagrangian relaxation formulation is depicted in Figure
\ref{lrformulation}. The Lagrangian dual is given by
\begin{equation*}
    LD \colon= \max_{0 \leq \kappa^C, 0   \leq \kappa^S } LR(\kappa^C,\kappa^S)
\end{equation*}
and this paper solves it using the Surrogate Lagrangian relaxation method proposed by \cite{bragin} to solve it.

\begin{figure}[!t]
\begin{align*}
& \min_{\textbf{x}, \textbf{y}} \quad & & \sum_{s \in S} \sigma_s*y_s+  & &\\
& & & \sum_{l \in L} \sum_{r \in R} \theta_l*m_l*v_r*x_{lr} + & & \\
& & & \sum_{l \in L} \kappa^C_l*(\sum_{r\in R}v_r*x_{lr} - c_l*y_{s_l}) + & & \\
& & &  \sum_{l \in L}\sum_{r\in R}  \kappa^S_{lr}*(x_{lr} - y_{s_l}) & & \\
&   \text{s.t. } & & \{\textbf{x}\} \in  \mathcal{G} \label{lag1} & & \\
 &  & & \textbf{x} \in \{0,1\}^{|L|X|R|}, \quad \textbf{y} \in \{0,1\}^{|S|} & &
\end{align*}
\caption{The Lagrangian Relaxation Formulation $LR(\kappa^C,\kappa^S)$.}
\label{lrformulation}
\end{figure}

%In each iteration of the sub-gradient algorithm the optimal solution of Lagrangian relaxation is found using the equations shown in Figure \ref{LRsolveeqn}. In the Lagrangian relaxation formulation the coupling constraints are relaxed, therefore the $\textbf{x}$ and $\textbf{y}$ variables are independent of each other. Among the $\textbf{x}$ variables the variables for each request $r$ i.e. $\textbf{x}^r$ are independent of each other. The  $\textbf{x}$ variables must satisfy the constraints $\mathcal{G}$, therefore the optimal value of $\textbf{x}^r$ variables can be found using $TDSPP$ algorithm, where the edge cost labels $\beta_{lr}$ = $\theta_l*m_l*v_r+\kappa^C_l*v_r+\kappa^S_{lr}$. The $\textbf{y}$ variables for each schedule $s$ i.e. $y_s$ are independent of each other, therefore the optimal value of $y_s$ variable can be found using the Equation (\ref{eq:ysoptimal}).
%\begin{figure}
%   \begin{subequations}\label{eq:lropt}
%\begin{equation}
%   \textbf{x}^r_{*} \gets  TDSP (G^*_r,\beta,\infty,1,0)  \label{lropt1} \\
%\end{equation}
%\begin{equation}
% y^{*}_s = \argmin_{y = \{0,1\}}\{ [\sigma_s - \sum_{l \in L^s} \sum_{r \in R}( \kappa^C_l *c_l + \kappa^S_{lr} ) ] * y  \} \\ \label{eq:ysoptimal}
%\end{equation}
%\end{subequations}
%\caption{Equations for optimal solution of $LR(\kappa^C,\kappa^S)$} \label{LRsolveeqn}
%\end{figure}

\section{Computational Results}
\label{section:experiments}

This section evaluates the performance of the algorithms and models
proposed in the paper. The evaluation considers two important settings:
\begin{enumerate}
\item \textbf{Tactical Planning}: Solve the TPOSSP for a full set of
  requests. During planning, many requests typically have base paths that
  serve as the initial basic feasible solution for the problem.
\item \textbf{Real-Time Operations}: Find routes for new requests on
  the day of operations while keeping the existing paths found during
  the TPOSSP planning step.
\end{enumerate}

\noindent
It is important to highlight three important facts. First, as mentioned
earlier, even after the proposed network reduction, the Arc-MIP still
contains hundreds of millions variables and is essentially
computationally intractable. Second, neither the planning optimization
nor the real-time optimization has been considered by the
industrial partner so far; they represent opportunities for
significant cost savings. In current practice, the TPOSSP is mostly a
manual process, where an existing base plan is incrementally updated
to accommodate new requests. Third, the scale of tactical planning is
substantially larger: it must handle all the requests of the network but has
more computational time. In contrast, real-time operations must solve
the problem quickly (e.g., within a minute) for relatively fewer requests.

\subsection{Instance Data}

The data used in the experiments is realistic and sourced from one of
the largest LSPs in the world.  The experiments includes 6 test cases
for planning (1.1 to 1.6) and 15 test cases for real-time operations
(2.1 to 2.15). A summary of the data is shown in Table
\ref{tab:data_summary}. Instance 1.1 covers half the network, while
instance 1.4 covers the full network. The real-time instances consist
of new requests and all the legs in the network with available
capacity. The real-time instances have wider time windows compared to
the planning instances. Four additional test instances for planning
have expanded time windows to measure the impact of increased time
flexibility on cost savings. Instances 1.2 and 1.3 use the same input
data as 1.1, with the window size expanded by 30 minutes in instance
1.2 (15 minutes before and 15 after) and by 1 hour in instance
1.3. Similarly, instances 1.5 and 1.6 use the same input data as 1.4,
with the window size expanded by 30 minutes and 1 hour
respectively. Table \ref{tab:data_summary} also provides the run time
of the EALSP algorithm to find the sub-network for each scenario.

\begin{table}[!t]
\small
\caption{Summary of Data\label{tab:data_summary}}
{\begin{tabular}{@{}>{\raggedright\arraybackslash}p{0.07\linewidth}@{\quad}>{\centering\arraybackslash}p{0.31\linewidth}@{\quad}>{\centering\arraybackslash}p{0.31\linewidth}@{}>{\raggedright\arraybackslash}p{0.2\linewidth}}
\hline
Instance& \# Requests x \# Schedule-legs& Average window size (hours)&EALSP run time (s)
 \\ 
\hline
 1.1&    16.7 Billion&14.4
&8000
 \\
1.2&    16.7 Billion&14.9&9000
 \\
 1.3&   16.7 Billion&15.4&
9600
 \\
1.4&  74 Billion & 9.0&20,580
 \\
1.5& 74 Billion & 9.5&22,980
 \\
1.6& 74 Billion & 10.0&21,000 \\
\hline
2.1&    2.34 Million &34.6&6
 \\
2.2&    11.7 Million &21.8&12
 \\
2.3&    15.7 Million &20.9&19
 \\
2.4&    36.6 Million &25.2&36
 \\
2.5&    200 Million &43.0&60
 \\
2.6&    127.5 Million &31.9&45
 \\
2.7&    47.8 Million &36.0&24
 \\
2.8&    63.7 Million &39.6&18
 \\
2.9&    66.6 Million &41.0&18
 \\
2.10&    33.9 Million &57.1&11
 \\
2.11&    56.8 Million &47.3&24
 \\
2.12&    46.6 Million &18.2&18
 \\
2.13&    50.5 Million &27.4&18
 \\
2.14&  26.4 Million & 40.1&54
 \\
2.15&  43.6 Million & 79.8&42
 \\
 \hline
\end{tabular}}{}
\end{table}

\subsection{Algorithm Settings}

The column-generation has four parameter values {\em Paths}, {\em
  NumIterations}, {\em MaxCost} and {\em Mode}. The {\em Paths}
parameter controls the number of paths generated per iteration and per
request. The {\em NumIterations} parameter controls the number of
iterations of the column generation. The {\em MaxCost} parameter
controls the maximum cost of the paths found by TDSPP algorithm, where
the cost of the paths range from the optimal cost (lowest possible) to
{\em maxCost}. The {\em Mode} parameter can either be {\em Standard}
or {\em Stabilized}, where the {\em Standard} mode runs the DW
decomposition for column generation, and the {\em Stabilized} mode
runs WDW for column generation. The parameter values used in the
experiments are given in Table \ref{tab:paramsettings}. The impact of
the {\em Paths} parameters was evaluated through experiments. Having
higher number of paths gives faster convergence (fewer iterations) at
the cost of increased run time per iteration. For large planning
datasets, the stabilization method has a significant impact on
run times and solution quality and the impact was evaluated
experimentally.  The lower bounds for real-time operations was
obtained by the exact Arc-MIP formulation. The lower bounds of the
planning instances were obtained by solving the Lagrangian dual, the
algorithm was ran until there was no improvement in objective for over
50 iterations. During planning, the CG-heuristic leverages the paths
assigned to existing requests in the base plan to seed the column
generation. By default, the column generation runs for 50
iterations. The experiments were conducted on a 24-core AMD EPYC 9334
processor and 786 GB of RAM, using the Gurobi 10.2 MIP solver.

\begin{table}[!t]
\small
\caption{Parameter Settings\label{tab:paramsettings}}
          {\begin{tabular}{@{}>{\raggedright\arraybackslash}p{0.4\linewidth}@{\quad}>{\centering\arraybackslash}p{0.25\linewidth}@{\quad}>{\centering\arraybackslash}p{0.25\linewidth}@{\quad}}
              \hline Parameter& Real-Time & Planning \\ \hline
              \textbf{CG-heuristic Parameters:} &\\ $Paths$& 50 &5,
              10, 20, 30, 40, 50 \\ $NumIterations$ &50 & 50
              \\ $MaxCost$ &0 &0 \\ $Mode$& Standard & Standard,
              Stabilized \\ \textbf{Gurobi Parameters for Arc-MIP:} &
              &\\ $MIPGap$ &0.05\% &\\ $TimeLimit$ &18,000 seconds
              &\\ \hline
\end{tabular}}{}
\end{table}

% The parameters used were $paths = 30$ and $numIterations = 20$ for scenarios 2.1 to 2.3, and $paths = 50$ and $numIterations = 20$ for scenarios 2.4 to 2.6.

\subsection{Business Benefits for Planning}

\begin{figure}[!t]
\begin{minipage}{0.48\textwidth}
\begin{tikzpicture}
  \begin{axis}[
    font=\footnotesize,   
    boxplot/draw direction=x,
    ytick = {1,2,3,4},
    yticklabels={Empty miles, Total miles, Schedule, Total cost},
    xlabel={\% Reduction from base},
    xticklabel={$\pgfmathprintnumber{\tick}\%$},
    height=8cm,
    width=8cm,
    ymajorgrids=true
  ]

    \addplot+[
      boxplot prepared={
        median=9.7,
        upper quartile=10.2,
        lower quartile=9,
        upper whisker=10.3,
        lower whisker=8.5
      },
      color=black!70
    ] coordinates {};  % Outlier
    \addplot+[
      boxplot prepared={
        median=3.2,
        upper quartile=3.5,
        lower quartile=3,
        upper whisker=3.7,
        lower whisker=2.8
      },
      color=black!70
    ] coordinates {};  % Outlier
 \addplot+[
      boxplot prepared={
        median=3.1,
        upper quartile=3.5,
        lower quartile=2.8,
        upper whisker=3.9,
        lower whisker=2.8
      },
      color=black!70
    ] coordinates {};
     \addplot+[
      boxplot prepared={
        median=2.6,
        upper quartile=2.9,
        lower quartile=2.3,
        upper whisker=3.2,
        lower whisker=2.3
      },
      color=black!70
    ] coordinates {};
  \end{axis}
\end{tikzpicture}
\end{minipage}
\hfill
\begin{minipage}{0.3\textwidth}
\begin{tikzpicture}
  \begin{axis}[
    font=\footnotesize,   
    boxplot/draw direction=y,
    xtick = {1},
    xticklabels={Optimality gap},
    ylabel={\% gap from LB},
    yticklabel={$\pgfmathprintnumber{\tick}\%$},
    height=8cm,
    width=4cm,
    ymajorgrids=true
  ]

    \addplot+[
      boxplot prepared={
        median=4.7,
        upper quartile=5.3,
        lower quartile=3.9,
        upper whisker=5.7,
        lower whisker=3.7
      },
      color=black!70
    ] coordinates {};
  \end{axis}
\end{tikzpicture}
\end{minipage}
\caption{Planning Results: the figure on the left shows the reduction
  in cost, total miles, schedules and empty miles compared to the base
  solution. The figure on the right is a box plot of the optimality
  gap achieved by the CG-heuristic for all the planning instances.}
\label{fig:scene2}
\end{figure}
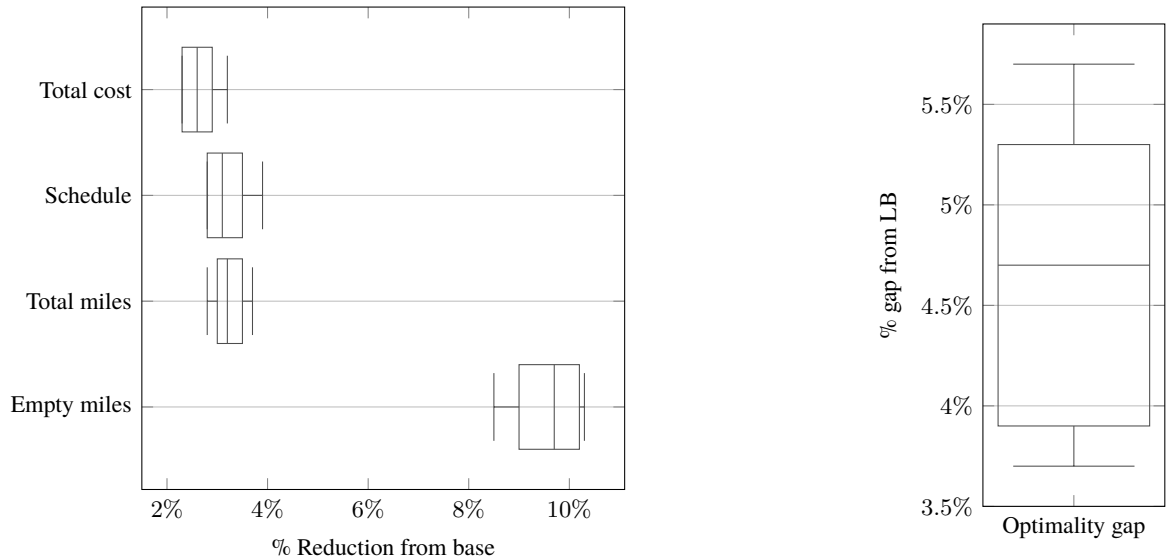

The best performance of CG-heuristic was obtained by using the optimal
sub-network and stabilized column generation. The cost reductions
obtained by CG-heuristic are calculated by taking the difference in
cost of the solutions proposed by CG-heuristic and the available base
solution.  Figure \ref{fig:scene2} depicts the results: cost
reductions range from 2.3\% to 3.2\%, {\em which translates into
  savings of tens of millions of dollars a year}.  Empty miles, which
are miles driven without any trailer, are reduced by 8.5\% to
10.3\%. The CG-heuristic achieves an optimality gap between 3.7\% and
5.7\%. The run times of CG-heuristic range from from 1.7 to 10.3
hours. Approximately 90\% of the total runtime is spent solving the
restricted master problem, 8\% solving the subproblem, and 2\% solving
the final MIP.

\subsection{The Benefits for Real-Time Operations}

\begin{figure}[!t]
     \noindent
     \begin{minipage}{0.48\textwidth}
\begin{tikzpicture}
  \begin{axis}[
   font=\footnotesize, 
    boxplot/draw direction=y,
    xticklabels={CG-heuristic MIP Gap},
    ylabel={\% Optimality Gap of CG-heuristic},
    height=8cm,
    width=8cm,
    ymajorgrids=true
  ]

    \addplot+[
      boxplot prepared={
        median=1.7,
        upper quartile=3,
        lower quartile=0.6,
        upper whisker=5.6,
        lower whisker=0
      },
      color=black!70
    ] coordinates {};  % Outlier

  \end{axis}
\end{tikzpicture}
\end{minipage}
\hfill
\begin{minipage}{0.48\textwidth}
\begin{tikzpicture}
  \begin{axis}[
   font=\footnotesize, 
    boxplot/draw direction=y,
   % ymode=log,                        % <--- This makes the y-axis logarithmic
    xtick={1,2},
    xticklabels={Arc-MIP, CG-heuristic},
    ylabel={Run time in seconds},
    height=8cm,
    width=8cm,
    ymajorgrids=true
  ]

    \addplot+[
      boxplot prepared={
        median=31,
        upper quartile=271,
        lower quartile=8.9,
        upper whisker=2369,
        lower whisker=1.8
      },
      color=black!70
    ] coordinates {};  % Outlier

    \addplot+[
      boxplot prepared={
        median=17,
        upper quartile=27,
        lower quartile=10.3,
        upper whisker=61,
        lower whisker=8.9
      },
      color=black!30
    ] coordinates {};        % No outliers

  \end{axis}
\end{tikzpicture}
\end{minipage}
\caption{Real-time results: the figure on the left is the optimality gap reached by CG-heuristic and the figure on the right is the run time of CG-heuristic and the time taken by Arc-MIP with the optimal sub-network to reach the same optimality gap}
\label{fig:scene1}
\end{figure}
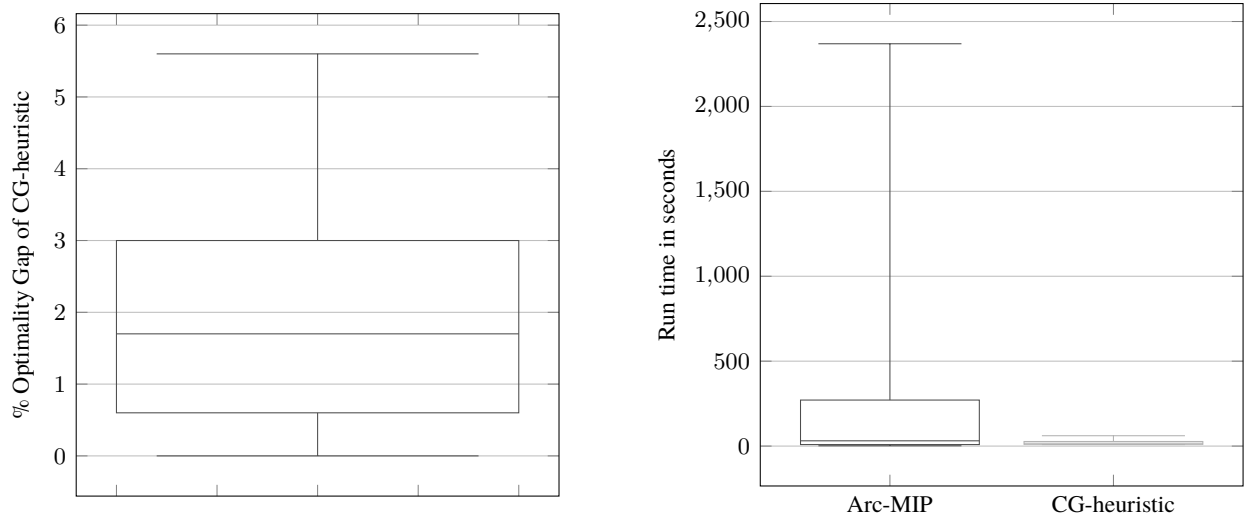

{\em One of the key outcomes of this paper is to demonstrate that the
  CG-heuristic can be used during real-time operations to schedule new
  requests.} In this setting, it is important to get high-quality
solutions quickly.  The computational results compare the runtime
performance of the CG-heuristic (run up to 50 iterations) and Arc-MIP
to obtain a solution within a given optimality gap. The left plot in
Figure \ref{fig:scene1} shows the optimality gap of the CG-heuristic
when compared to the lower bound found by solving the Arc-MIP model
(running for up to 5 hours). The right plot in Figure \ref{fig:scene1}
contrasts the run times of the CG-heuristic with those of Arc-MIP which
is terminated when it reaches the optimality gap achieved by the
CG-heuristic. The CG-heuristic provides a solution within 3\% of the
true optimal solution in 75\% of test cases. Its run time is on average
85\% faster than Arc-MIP. It is almost always faster, except in a few
small instances. Overall, these results indicate that the CG-heuristic
is now mature for real-time operations.

\subsection{The Benefits of the Optimal Sub-Network}

The small real-time instances 2.1 to 2.4 were chosen to demonstrate
the impact of computing the optimal sub-network as a pre-processing
step. It was found that the run times decreased by 83\%-99\% when
solving Arc-MIP on the optimal sub-network.  The run times of the
CG-heuristic decreased by 53\%-82\% when using the optimal
sub-network.

For planning, finding the optimal sub-network takes about 5 hours for
the full-network, as shown in Table \ref{tab:data_summary}. However,
the sub-network algorithm does not need to be run every time a
planning step comes. Indeed, {\em the optimal sub-network has already
  been computed for existing requests and it suffices to update the
  optimal sub-network whenever there is a change in a schedule-leg or
  a new request is added.} In practice, when a planning solution is
requested, the sub-network algorithm runs only for the newly added
requests, as the optimal sub-network for pre-existing requests will
already be available.

\subsection{The Impact of Stabilization and Varying Paths}

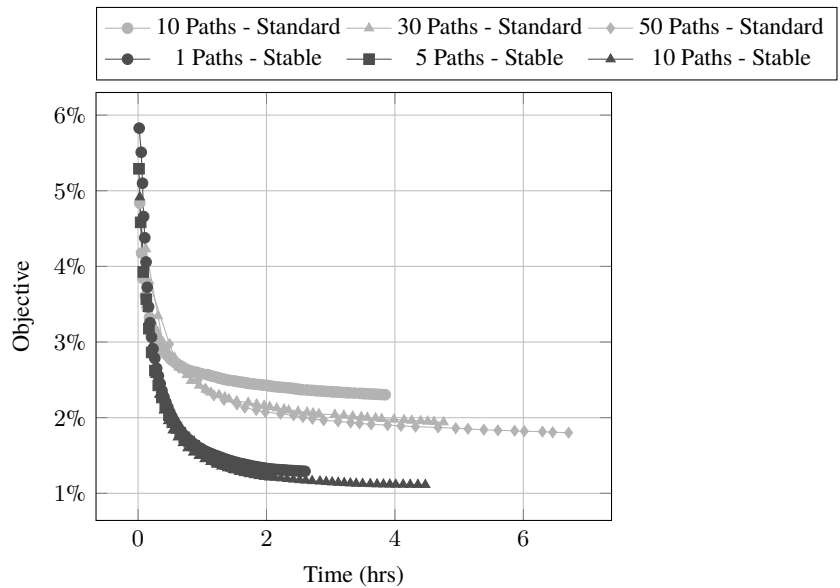
\begin{figure}[!t]
\definecolor{ForestGreen}{RGB}{34,139,34}
\definecolor{DarkOrange}{RGB}{255,140,0}
  \noindent
  \begin{tikzpicture}
  \begin{axis}[
   font=\footnotesize, 
    xlabel={Time (hrs)},
    ylabel={Objective},
    yticklabel=\pgfmathprintnumber{\tick}\%% <-- Add % symbol
    ,
    grid=both,
    legend style={at={(0,1.2)}, anchor=north west, legend columns=3}
  ]
\addplot[
  color=black!30,
  mark=*,
] table[
  x=time, y={10},
  col sep=comma
] {info/path_experiment.csv};
\addlegendentry{10 Paths - Standard}

\addplot[
  color=black!30,
  mark=triangle*,
] table[
  x=time, y={30},
  col sep=comma
] {info/path_experiment.csv};
\addlegendentry{30 Paths - Standard}

\addplot[
  color=black!30,
  mark=diamond*,
] table[
  x=time, y={50},
  col sep=comma
] {info/path_experiment.csv};
\addlegendentry{50 Paths - Standard}

\addplot[
  color=black!70,
  mark=*,
] table[
  x=time, y={1_stable},
  col sep=comma
] {info/stable_path_experiment.csv};
\addlegendentry{1 Paths - Stable}

\addplot[
  color=black!70,
  mark=square*,
] table[
  x=time, y={5_stable},
  col sep=comma
] {info/stable_path_experiment.csv};
\addlegendentry{5 Paths - Stable}

\addplot[
  color=black!70,
  mark=triangle*,
] table[
  x=time, y={10_stable},
  col sep=comma
] {info/stable_path_experiment.csv};
\addlegendentry{10 Paths - Stable}

  \end{axis}
\end{tikzpicture}

\caption{The Impact of stabilization and varying paths generated per request per iteration of column generation}
\label{stableandpaths}
\end{figure}

The $Paths$ parameter specifies the maximum number of paths to
generate per request per iteration of column generation. The impact of
changing $Paths$ is measured on standard column generation. Figure
\ref{stableandpaths} shows a reduction in the objective over time when
increasing the $Paths$. With increasing number of paths, the
performance of standard column generation improved significantly,
whereas the stabilized column generation did not see a significant
improvement in performance.  The stabilized column generation performs
strictly better than standard column generation with 50 paths; it
converges quicker and flattened out on average with an objective that
is 1\% lower.

\section{Conclusion}
\label{section:conclusion}

This paper addresses the trailer path optimization with schedule
services (TPOSSP) problem, a common challenge faced by LSPs. A review
of the literature revealed that the TPOSSP has not been
explored in the literature. However, mathematical models with
structure similar to TPOSSP have been solved for rail network
optimization and SND problems in the literature.

Due to the large-scale nature of practical instances, it was
imperative to reduce the size of the networks considered by the
solver. The paper introduced the concept of an optimal sub-network to
prune out legs that never lie on any feasible path of a request.  This
pruning step reduces the run times by 50\% - 82\%. However, this step
alone was not sufficient for the tactical planning/re-planning
setting.  Indeed, the MIP models still contain hundreds of millions
variables. To address this complexity, the paper proposed a stabilized
Dantzig-Wolfe decomposition whose pricing problem is a time-dependent
shortest path algorithm.

The CG-heuristic was tested on real-life datasets. It was found that
for real-time operations, the heuristic gave results within 3\%
optimality gap for 75\% of cases and within 5.3\% gap for all the
cases (Figure \ref{fig:scene1}) and runs within 1 minute. For
planning/re-planning purposes, the CG-heuristic give an optimality gap
below 5.3\% for 75\% of cases and run times ranging from 1.7 to 10.3
hours.  It also reduces operation costs by 2.3\% - 3.2\% and empty
miles by 8.5\% to 10.3\% (\ref{fig:scene2}).

In summary, this paper makes two fundamental contributions: it shows
that (1) the CG-heuristic can solve planning/re-planning of the scale
encountered in practice, producing high-quality solutions to problems
involving hundreds of millions variables, even after valid network
pruning; 2) the CG-heuristic is fast enough to be deployed in
real-time settings where it can optimize a set of additional requests
with a minute; and (3) the savings resulting from the CG-heuristic are
in the tens of millions of dollars a year.

% Appendix here
% Options are (1) APPENDIX (with or without general title) or
%             (2) APPENDICES (if it has more than one unrelated sections)
% Outcomment the appropriate case if necessary
%
% \begin{APPENDIX}{<Title of the Appendix>}
% \end{APPENDIX}
%
%   or
%
% \begin{APPENDICES}
% \section{<Title of Section A>}
% \section{<Title of Section B>}
% etc
% \end{APPENDICES}

% Acknowledgments here
\newpage

\appendix
\section{Proofs}

\begin{proof} {\textbf{Proof of Lemma \ref{feaspathlemma}:}} 
    A path can be constructed passing through $l$ in two phases. It must be true that there is a feasible path from origin of request $r$ to $i$, because if not then by definition $eat^r_i = \infty$ and  $eat^r_i \leq lst^r_i$ will not hold. Similarly, there is a feasible path from $j$ to destination of request $r$, because if not then by definition $lst^r_i = -\infty$ and  $eat^r_i \leq lst^r_i$ will not hold. Let the path from origin of request $r$ to $i$ be given by $P^i$ and path from hub $j$ to destination of request be given by $P^j$. Connect $P^i$ and $P^j$ by leg $l$. $P^i$ $\rightarrow$ $l$ $\rightarrow$ $P^j$ is a feasible path, because it is true that $l$ departs at or after $eat^r_i$ which is the arrival time of $P^i$ and $l$ arrives at $j$ at or before $lst^r_j$ which is the departure time of $P^j$.
\end{proof}
\begin{proof} {\textbf{Proof of Theorem \ref{thm1}:}} 
Consider $E^r = \Phi$ then by definition of $eat^r$ and $lst^r$ it is clear that there is no feasible path, so therefore $E^r_{*} = \Phi$. This is true because, if there is a feasible path it must be true that $eat^r \leq lst^r$ for all the hubs of the path which would imply $E^r \neq \Phi$. 

Consider $E^r \neq \Phi$, first  $E^r_{*} \subseteq E^r$ is proven and after that  $E^r \subseteq E^r_{*}$ is proven.
\begin{itemize}
    \item  \underline{ $E^r_{*} \subseteq E^r$} : Let schedule-leg $l \in E^r_{*}$, with origin hub $i$ and destination hub $j$. It is true from the definition of $eat^r$ and $lst^r$ that for every $h \in H^r_{*}$, $eat^r_h \leq lat^r_{h}$. Therefore it must be true that  $eat^r_i \leq t^s_{l}$ and $t^e_{l} \leq lst^r_j$.
    \item \underline{ $E^r \subseteq E^r_{*}$} : Let leg $l \in E^r$, it will be shown that $l$ lies on a path from origin to destination of request $r$. All legs in $l \in E^r$ satisfy the condition $eat^r_i \leq t^s_l \leq lst^r_j$, where $i$ is the origin of leg $l$ and $j$ is the destination of leg $l$. Applying Lemma \ref{feaspathlemma} to every $l \in E^r$ it is true that there is at least one feasible path passing every leg in $E^r$. By the definition of optimal sub-network it is true that $E^r_{*}$ contains all the legs used by a feasible path, therefore it must be true that $E^r \subseteq E^r_{*}$.
\end{itemize}
Therefore $E^r_{*}=E^r$.
\end{proof}

\begin{proof} {\textbf{Proof of Lemma \ref{lemma2}:}} 
The lemma is proven by contradiction. Let's suppose that there exists a leg $l$ with departure time $t^s_l$ that is not critical but lies on an all the optimal paths. Let the next critical point be $t^{'}$ and the corresponding leg be $l^{'}$. Formally, $t^{'} = \inf_{t \in criticalpoints} t|t\geq t^s_l$. Where $criticalpoints$ is the set of critical points of function $f^r_{d_l}$. It is known that $f^r_l(t)$ is step-wise left continuous non-decreasing function, therefore it is true that $f^r_{d_l}(t^s_l) = f^r_{d_l}(t^{'})$. Any optimal path that passes through $l$ must arrive at $d_l$ at or before $t^s_l$. But, since $f^r_{d_l}(t^s_l) = f^r_{d_l}(t^{'})$ because of step-wise left-continuous property it means that an optimal path can alternately use leg $l^{'}$ which is a critical leg instead of leg $l$, which contradicts the fact that $l$ lies on all optimal paths. Therefore, it is true that legs departing on non-critical points can be removed without affecting the optimal solution.
\end{proof}

\clearpage
% References here (outcomment the appropriate case)
%
% CASE 1: BiBTeX used to constantly update the references
%   (while the paper is being written).
%\bibliographystyle{informs2014trsc} % outcomment this and next line in Case 1
%\bibliography{<your bib file(s)>} % if more than one, comma separated

\bibliographystyle{plainnat}  % or unsrt, abbrv, etc.
\bibliography{TPOSSP} % if more than one, comma separated

% CASE 2: BiBTeX used to generate mypaper.bbl (to be further fine tuned)
%\input{mypaper.bbl} % outcomment this line in Case 2

%If you don't use BiBTex, you can manually itemize references as shown below.

%\bibliographystyle{nonumber}

%%%%%%%%%%%%%%%%%
\end{document}